\documentclass[final,3p,11pt, times]{article}
\usepackage{authblk}
\usepackage{amsmath}
\usepackage{amssymb,latexsym}
\usepackage[german,english]{babel}
\usepackage{url}
\usepackage{graphicx}
\usepackage{gastex}
\usepackage{longtable}
\usepackage{lscape}
\usepackage{tabularx}
\usepackage{multicol}
\usepackage{verbatim}
\usepackage{multirow}
\usepackage{graphicx}
\usepackage{epsfig,graphics,graphicx}
\usepackage{geometry}
\usepackage[usenames]{xcolor}
\usepackage{tikz}
\hbadness=\maxdimen

\newtheorem{lem}{{\bf Lemma}}[subsection]
\newtheorem{thm}{{\bf Theorem}}[section]
\newtheorem{prop}{{\bf Proposition}}[section]

\usepackage{chngcntr}
\counterwithout{equation}{section}
\begin{document}
	
	\title{On the $A_{\alpha}$ and $RD_{\alpha}$ matrices over certain groups}
	
	\author[1]{Yogendra Singh}
	\author[1]{Anand Kumar Tiwari}
	\author[2]{Fawad Ali}
	
	\affil[1]{\small Department of Applied Science, Indian Institute of Information Technology, Allahabad 211015, India}

	\affil[2]{Institute of Numerical Sciences, Kohat University of Science and Technology, Kohat 26000, KPK, Pakistan}
	
	\maketitle
	
	\hrule
	
\begin{abstract}
The power graph $G = P(\Omega)$ of a finite group  $\Omega$ is a graph with the vertex set $\Omega$ and two vertices $u, v \in \Omega$ form an edge if and only if one is an integral power of the other. Let $D(G)$, $A(G)$, $RT(G)$, and $RD(G)$ denote the degree diagonal matrix, adjacency matrix, the diagonal matrix of the vertex reciprocal transmission, and Harary matrix of the power graph $G$ respectively. Then the $A_{\alpha}$ and $RD_{\alpha}$ matrices of $G$ are defined as $A_{\alpha}(G) = \alpha D(G) + (1-\alpha)A(G)$ and $RD_{\alpha}(G) = \alpha RT(G) + (1-\alpha)RD(G)$. In this article, we determine the eigenvalues of $A_{\alpha}$ and $RD_{\alpha}$ matrices of the power graph of group $ \mathcal{G} = \langle s,r \, : r^{2^kp} = s^2 = e,~ srs^{-1} = r^{2^{k-1}p-1}\rangle$. In addition, we calculate its distant and detotar distance degree sequences, metric dimension, and strong metric dimension.  
\end{abstract}

\smallskip



\hrule
\section{Introduction}\label{s1}
The study of various algebraic structures through several associated graphs has been a fascinating area of research in the last few years, see \cite{ ah(2009), ain(2017), ash(2016), cg(2011), dk(2021)}. In the present study, we focus on the associated power graphs of group $\mathcal{G}$. The notion of power graph (directed) was introduced by Kelarev and Quinn \cite{kq(2002)} for semigroups as follows: if $S$ is a semi-group then the power graph $P(S)$ of $S$ is a graph that has $S$ as the vertex and $uv$ is an edge in $P(S)$ from a vertex $u$ to another vertex $v$ if $u = v^n$ for some $n \in \mathbb{N}$. Motivated by this, Chakrabarty et al. \cite{cgs(2009)} defined an undirected power graph for a group in which two distinct vertices are adjacent if and only if  one can be written as an integral power of the other. Further, the authors proved that the power graph $P(\Omega)$ of a group $\Omega$ is complete if and only if the order of $\Omega$ is $1$ or $p^n$ for some prime $p$ and $n \in \mathbb{N}$. For more results, we refer to survey articles \cite{kscc(2021)} and \cite{akc(2013)}.

By a graph $G = G(V, E)$, we mean an undirected simple graph with the vertex set $V$ and edge set $E$. The cardinality of the set $V$ is called the order of the graph $G$ and is denoted as $|V|$. The sum of two connected graphs $G_1$ and $G_2$, denoted as $G_1 + G_2$, is a graph with the vertex set $V_1 \cup V_2$ and edge set $E_{1} \cup E_{2} \cup \{ u_i \sim v_i : u_i \in V_{1},~ v_i \in V_{2}\}$, respectively. We denote the disjoint union of $n$ copies of $G$ by $nG$ and complete graph on $n$ vertices by $K_{n}.$ For $v_i, v_j \in V$, $v_i \sim v_j$ means $v_i$ and $v_j$ are adjacent, otherwise $v_i \nsim v_j$.
For $u_i \in V$, the sets $N(u_i) = \{ v_i \in G \ | \ u_i \sim v_i\}$ and $N[u_i] = N(u_i) \cup \{u_i\}$ are called open and closed neighborhood of $u_i$ respectively. The degree $\deg(v_i)$ of a vertex $v_i$ is $|N(v_i)|.$ The vertices $u_i$ and $v_i$ are called twins if $N[u_i] = N[v_i]$ or $N(u_i) = N(v_i).$ A subset $T$ of $V$ is called a twin-set if every distinct pair of vertices in $T$ are twins.

A subset $U$ of $G$ is called a vertex cover of $G$ if every edge of $G$ is incident with at least one vertex of $U.$ The smallest cardinality of a vertex cover of $G$ is called the vertex covering number of $G$ and is denoted by $\beta(G)$.

For an ordered subset $R=\{r_{1}, r_{2}, \ldots , r_{k}\}$ of $G$ and $u_i \in G$, the vector $r(u_i|R)=(d(u_i,~ r_{1}),~ d(u_i,~ r_{2}), \ldots ,~ d(u_i, ~r_{k}))$ denotes the representation of $u_i$ as the $k$-vector with respect to $R$. The set $R$ is called a resolving (or locating) set of $G$ if for any two distinct vertices $u_i,v_i \in G$, $r(u_i|R) \neq r(v_i|R)$. The minimal resolving set of $G$ is called the basis of $G$ and its cardinality is the metric dimension of $G$, represented by $\psi(\Gamma)$. Slater explored the ideas of metric dimension, which were then researched separately by several researchers, such as \cite{ hm(1976), cej(2020)}.

A pair of vertices $v_i,w_i \in G$ is said to be strongly resolve by a vertex $u_i \in G$ if there exists a shortest $v_i$-$u_i$ path containing $w_i$ or a shortest $w_i$-$u_i$ path containing $v_i.$ A subset $S$ of $G$ is called a strong resolving set if every pair of vertices of $G$ is strongly resolved by some vertex of $S.$ The smallest cardinality of a strong resolving set for $G$ is called the strong dimension of $G$ and is denoted by $sdim(G)$.

A vertex $u_i \in G$ is called maximally distant from $v_i \in G$ (denoted by $u_i \ MD \ v_i$) if for every $x \in N(u_i)$, $d(v_i,x) \leq d(u_i, v_i)$. We say that $u_i$ and $v_i$ are mutually maximally distant (denoted by $u_i \ MMD \ v_i$) if $u_i \ MD \ v_i$ and $v_i \ MD \ u_i$. The strong resolving graph $G_{SR}$ of a connected graph $G$ ia graph having vertex set $V_{SR} = V$ and $u_iv_i \in E_{SR}$ iff $u_i \ MMD \ v_i$.

The distance (detour distance) $d(u_i,v_i)$ ($d_D(u_i,v_i)$) between the vertices $u_i,v_i \in G$ is the length of a shortest (largest) $u_i$-$v_i$ path in $G$. The eccentricity (detour eccentricity) of a vertex $u_i \in G$ is the largest distance (detour distance) between $u_i$ and any other vertex of $G$, denoted as $ec(u_i) \, (ec_{D}(u_i))$. The diameter (detour diameter) of $G$, denoted by  $dia(G)\, (dia_D(G))$, is the largest eccentricity among all the vertices of $G$. Also, the eccentricity (detour eccentricity), denoted by $rad(G) \,(rad_D(G))$, is the smallest eccentricity.

Let $deg_j(u_i)(Deg_{j}(u_i))$ be the set of vertices at distance (detour distance) $j$ from $u_i \in  V(G)$. Then
$(deg_{0}(u_i)$, $deg_{1}(u_i)$, $deg_{2}(u_i),\ldots, deg_{ec(u_i)}(u_i)((Deg_{0}(u_i)$, $Deg_{1}(u_i)$, $Deg_{2}(u_i), \ldots, Deg_{ec_{D}(u_i)}(u_i))$ is called the distance (detour distance) degree sequence and is denoted by $dds(u_i) (dds_D(u_i))$. 
The $dds(G)\, (dds_D(G))$ of $G$ is the set of $dds(u_i)\, (dds_D(u_i))$, where $u_i \in V(G)$.

Let $A(G)$ and $D(G) = diag ( \deg(v_1), \deg(v_2), \ldots, \deg(v_{|V|}))$ denote the adjacency matrix and  degree diagonal matrix of a graph $G$ of order $n$ respectively.  Then the Laplacian matrix $L(G)$ and signless Laplacian matrix $Q(G)$ are defined as $L(G) = A(G) - D(G)$ and $Q(G) = A(G)+D(G)$ respectively.  Since the matrices $A(G), L(G)$, and $Q(G)$ are symmetric, so all of there eigenvalues are real. Hence the eigenvalues $ \lambda_i^A$'s, $\lambda_i^L$'s and $\lambda_j^Q$'s of $A(G), L(G)$ and $Q(G)$ can be written as $\lambda_1^A \geq \lambda_2^A \geq \ldots \geq \lambda_n^A$, $\lambda_1^L \geq \lambda_2^L \geq \ldots \geq \lambda_n^L$ and $\lambda_1^Q \geq \lambda_2^Q \geq \ldots \geq \lambda_n^Q$ respectively. We denote the multiplicity of an eigenvalue $\lambda$ by $m(\lambda).$ The largest eigenvalue $\lambda_1^A$ is called the spectral radius or spectral norm of $G.$

Nikiforov \cite{n(2017)} considered the matrix $A_{\alpha}(G)$ as  $A_{\alpha} = \alpha A(G) + (1 - \alpha) D(G)$, where $\alpha \in [0, 1]$. Note that $A_0(G) = A(G)$, $A_1(G) = D(G)$, $A_{\frac{1}{2}}(G) = \frac{1}{2} Q(G)$, and $A_{\alpha}(G) - A_{\beta}(G) = (\beta - \alpha) L(G)$. Thus the $A_{\alpha}$ matrix combines the spectral theories of the matrices $A(G)$, $D(G)$, $L(G)$, and $Q(G)$ including many other combinations of $\alpha$. For the recent progress on the $A_{\alpha}$ matrix, we refer to see \cite{ls(2021), lw(2022), lz(2022), fw(2022)}. 

Plavs$\acute{i}$c \cite{pntm(1993)} introduced the notion of Harary matrix (also called reciprocal distance matrix) of a graph $G$, as 
\begin{equation*}
RD(G) =(a_{ij})=\begin{cases}
          \frac{1}{d(u_i, u_j)} \quad &\text{if} \, u_i \neq u_j \\
          0 \quad &\text{otherwise.} \\
     \end{cases}
\end{equation*}
The reciprocal transmission $RT(v)$ of a vertex $v$ is defined as $RT(v) = \sum\limits_{u \in V(G) \setminus \{v\}} \frac{1}{d(v,u)}$. Then the vertex reciprocal transmissions matrix $RT(G)$ of $G$ is the diagonal matrix $RT(G) = dia(RT(v_1), RT(v_2), \ldots, RT(V_{|V|}))$. The matrices $RL(G) = RT(G)-RD(G)$ and $RQ(G) = RT(G) + RD(G)$ are called reciprocal distance Laplacian matrix and reciprocal distance signless Laplacian matrix. Motivated by the $A_{\alpha}(G)$, discussed above, Tian et al. \cite{tcc} constructed the matrix $RD_{\alpha}(G)$ as $RD_{\alpha}(G) = \alpha RT(G) + (1 - \alpha) RD(G)$, where $\alpha \in [0, 1]$, that combines the spectral properties of $RD(G)$, $RQ(G)$, $RT(G)$.  

Graph spectra have several applications in Quantum chemistry and in several mathematical branches such as group theory, Lie algebra, and ring theory. Over the last few years many researchers have shown their interest in calculating the spectrum of algebraic graphs, see \cite{ Rather2, Ali(2020), AL(2019), pgab(2020), 6, 6a, 15aa}.

Motivated by the above work, here, we attempt to study the spectra of power graphs of group $\mathcal{G}$ through the theory of $A_{\alpha}$ matrix and $RD_{\alpha}$ matrices. We proceed in the following manner. In Section 2, we describe the power graph of group $\mathcal{G}$. In Section 3, we determine the $A_{\alpha}$ eigenvalues of the power graph of group $\mathcal{G}$.  In Section 4, we determine the $SD_{\alpha}$ eigenvalues of the power graph of group $\mathcal{G}$.  In Section 5, we determine the power graph distant and detotar distance degree sequences, metric dimension, and strong metric dimension.

\section{Power graph of group $\mathcal{G}$}\label{s2}

In this article, we focus on the group $\Omega_{2,2^kp,2^{k-1}p-1}$ with the following representation: $$\Omega_{2,2^kp,2^{k-1}p-1} = \langle s,r \, : r^{2^kp} = s^2 = e,~ srs^{-1} = r^{2^{k-1}p-1}\rangle.$$

\smallskip

Here, we denote $\Omega_{2,2^kp,2^{k-1}p-1}$ by $\mathcal{G}$. Observe that $\mathcal{G}$ is a non-abelian group of order $2^{k+1}p)$, where  $p > 2$ is a prime number and $k \geq 2$. 

In Section \ref{s5}, we consider the following partition of the group $\mathcal{G} = \{H_0 , H_1, H_2, \linebreak H_3 \}$, where $H_0 = \{e, u = r^{2^{k-1}p}\},~
H_1 =  \{r^i : 1 \leq i \leq 2^kp-1\} \ \setminus \ \{u\},~
H_2 =  \{sr^{2t} {: 1 \leq t \leq 2^{k-1}p}\},$ and $H_3 = \{sr^{2j+1} : 0 \leq j \leq 2^{k-1}p-1\}$.

For the above structure, we have: 

\begin{thm} For the group $\mathcal{G}$, we have
$$P(\mathcal{G}) = P \left(\mathbb{Z}_{2^kp} \right) \cup 2^{k-1}p K_2 \cup 2^{k-2}p \left((P \langle r^{2^{k-1}p} \rangle) + K_2 \right).$$
\end{thm}

\section{$A_{\alpha}$ eigenvalues of the power graph $P(\mathcal{G})$}\label{s3}

In this section, we determine the $A_{\alpha}$ the power graph $P(\mathcal{G})$.

\smallskip
For a column vector $X = (x_1, x_2, 
\ldots, x_n)^T \in \mathbb{R}^n$, we have

$$\langle A_{\alpha}X, X\rangle = (2\alpha-1)\sum\limits_{u \in V}^{} x^2_{u}deg(u) + (1-\alpha) \sum\limits_{uv \in E}^{} (x_u + x_v)^2,$$

and if $\lambda$ is an $A_{\alpha}$ eigenvalue of a graph $G$ corresponding to the vector $X (\neq 0)$, then

\begin{equation} \label{e6.1}
\lambda x_i = \alpha \cdot deg(u_i)x_i + (1-\alpha)\sum_{u_{i}u_{j} \in E(G)}^{} x_j.
\end{equation}

Now, recall some results which are useful to determine some  $A_{\alpha}$ eigenvalues of $G$.

\begin{prop} \cite{bofl(2019)}  \label{p1} Let $u_i$ and $u_{j_t}$, $1 \leq t \leq l$, are twin vertices in a graph $G$ on $n \geq 2$ vertices.
	
	\begin{enumerate}
		\item If $u_i \nsim u_{j_t}$, then $\alpha \cdot deg(u_{i})$ is an  $A_{\alpha}$ eigenvalue of $G$ and $m(\alpha \cdot deg(u_{i})) \geq l.$
		
		\item If $u_i \sim u_{j_t}$, then $\alpha (deg(u_{i})+1)-1$ is an $A_{\alpha}$ eigenvalue of $G$ and $m(\alpha (deg(u_{i})+1)-1) \geq l.$
	\end{enumerate}

\end{prop}

\begin{thm} [\cite{BH}] \label{t6.1}
	If $\mathcal{A}, \mathcal{B}, \mathcal{C}$, and $\mathcal{D}$ are square matrices and $\mathcal{D}$ is invertible. Then \[\begin{vmatrix}
	\mathcal{A} & \mathcal{B}  \\
	\mathcal{C} & \mathcal{D} \end{vmatrix} = \begin{vmatrix}\mathcal{D}	\end{vmatrix} \begin{vmatrix} \mathcal{A} - \mathcal{B}\mathcal{D}^{-1}\mathcal{C} \end{vmatrix}. \]
\end{thm}

\begin{thm}
 The $A_{\alpha}$ matrix of $P(\mathcal{G})$ have the eigenvalues $\alpha$ with $m(\alpha)=2^{k-1}p-1$,  $\alpha \cdot 2^kp -1$ with $m(\alpha(2^kp -1)=2^{k}p-3$, $4\alpha -1$ with $m(4\alpha -1)=2^{k-2}p$, and $2\alpha +1$ with $m(2\alpha +1)=2^{k-2}p-1$ respectively, and the roots of the five degree polynomial $P(x)$ given in the proof. 
\end{thm}

\noindent{\textbf{Proof.}} The vertex set of $P(\mathcal{G})$ is
$\{1, r, r^2, \ldots, r^{2^kp-1}, s, sr^2, \ldots, sr^{2^kp-2}, sr, sr^3, \ldots, sr^{2^kp-1}\}$. Let $T_1 = \{s, sr^2, \ldots, sr^{2^kp-2}\}$, $T_2 = \{r, r^2, \ldots, r^{2^kp-1}\}\setminus \{r^{2^{k-1}p}\}$, $T_j = \{ rs^{2j+1}, rs^{2j+1+2^{k-1}p}\}$ for  $0 \leq j \leq 2^{k-2}p-1.$ Since the elements of $T_1$ are non adjacent twins and they have degree $1$, by Proposition \ref{p1}, we get $\alpha$ is an $A_{\alpha}$ eigenvalue with $m(\alpha)=2^{k-1}p-1.$ As the elements of $T_2$ are adjacent twins and they have degree $2^{k}p-1$, by Proposition \ref{p1}, we get $\alpha 2^kp -1$ is an $A_{\alpha}$ eigenvalue with $m(\alpha2^kp -1)=2^{k}p-3.$ Also, for a fixed $j$ the elements of $T_j$ are adjacent twins and they have degree 3, so by Proposition \ref{p1}, we get $4\alpha -1$ is an $A_{\alpha}$ eigenvalue with $m(4\alpha -1)=1.$ Since there are $2^{k-2}p$ sets of type $T_j$, so $4\alpha -1$ is an $A_{\alpha}$ eigenvalue with $m(4\alpha -1)=2^{k-2}p.$ Thus, we are remains with $2^{k-2}p+4$ eigenvalues. 

Let $X$ be the eigenvector of $A_{\alpha}$ such that $x_i = X(v_i)$ for $1 \leq i \leq 2^{k+1}p$. Then, clearly the components of $X$ that corresponds to $u \in T_1, 1, v \in T_2 \text{ and } r^{2^{k-1}p}$ are $x_1, x_2, x_3$ and $x_4$ respectively, the components of $X$ that corresponds to both $rs^{2j+1}$ and $rs^{2j+1+2^{k-1}p}$ is
$x_{j+5}$, where $0 \leq j \leq 2^{k-2}p-1$. Thus from Equation \ref{e6.1}, we have:
\begin{align*}
\lambda x_1 &= \alpha x_1 + (1-\alpha)x_2\\
\lambda x_2 &= 2^{k-1}p(1-\alpha) x_1 + (2^{k+1}p-1) \alpha x_2 +(2^kp-2)(1-\alpha)x_3 + (1 - \alpha)x_4 \\
&+ 2(1-\alpha)\sum_{j=0}^{j=2^{k-2}p-1}x_{5+j}\\
\lambda x_3 &= (1-\alpha) x_2 + (2^{k}p-3)(1-\alpha) x_3 +(2^kp-1)\alpha x_3 + (1 - \alpha)x_4\\
\lambda x_4 &= (1-\alpha) x_2 + (2^{k}p-2)(1-\alpha) x_3 +(3 \cdot 2^{k-1}p-1)\alpha x_4 + (1 - \alpha)x_4 \\
&+2(1-\alpha)\sum_{j=0}^{j=2^{k-2}p-1}x_{5+j}\\
\lambda x_5 &= (1-\alpha) x_2 + (1-\alpha) x_4 +3\alpha x_5 + (1 - \alpha)x_5\\
\lambda x_6 &= (1-\alpha) x_2 + (1-\alpha) x_4 +3\alpha x_6 + (1 - \alpha)x_6\\
&~~\vdots\\
\lambda x_{5+2^{k-2}p-1} &= (1-\alpha) x_2 + (1-\alpha) x_4 +3\alpha x_{5+2^{k-2}p-1} + (1 - \alpha)x_{5+2^{k-2}p-1}
\end{align*}

In the above system of equation, the coefficient matrix corresponding to the right side is 
 $$\mathcal{Y} = \begin{bmatrix}
\mathcal{A}_{4 \times 4} & \mathcal{B}_{4 \times 2^{k-2}p}  \\
\mathcal{C}_{2^{k-2}p \times 4} & \mathcal{D}_{2^{k-2}p \times 2^{k-2}p} \end{bmatrix}, \text{where}$$

$$\mathcal{A}_{4 \times 4} = \begin{array}{@{} c @{}}
	\begin{bmatrix}
	\begin{array}{ *{8}{c} }	
	\alpha  & 1-\alpha & 0 &  0 \\
	2^{k-1}p(1-\alpha) & (2^{k+1}p-1)\alpha & (2^{k}p-2)(1-\alpha) & 1-\alpha  \\
	0 & 1-\alpha & 2^{k}p+2\alpha-3 & 1-\alpha \\
	0 & 1-\alpha & (2^{k}p-2)(1-\alpha) & (3 \cdot 2^{k-1}p-1)\alpha \\
\end{array}
\end{bmatrix}\mathstrut
\end{array},$$
$$ \mathcal{B}_{4 \times 2^{k-2}p} = \begin{array}{@{} c @{}}
\begin{bmatrix}
\begin{array}{ *{8}{c} }
0 & 0 & \ldots & 0 \\
2(1-\alpha) & 2(1-\alpha) & \ldots & 2(1-\alpha)  \\
0 & 0 & \ldots & 0 \\
2(1-\alpha) & 2(1-\alpha) & \ldots & 2(1-\alpha)  \\
\end{array}
\end{bmatrix}, \ \mathcal{C}_{2^{k-2}p \times 4} = 
\begin{bmatrix}
\begin{array}{ *{8}{c} }
0 & 1-\alpha & 0 & 1-\alpha \\
0 & 1-\alpha & 0 & 1-\alpha \\
\vdots & \vdots  & \vdots & \vdots  \\
0 & 1-\alpha & 0 & 1-\alpha  \\
\end{array}
\end{bmatrix}\mathstrut
\end{array},$$ 

$$ \mathcal{D}_{2^{k-2}p \times 2^{k-2}p} = \begin{array}{@{} c @{}}
\begin{bmatrix}
\begin{array}{ *{8}{c} }
2\alpha +1 & 0 & \cdots & 0 \\
0 & 2\alpha +1 & \cdots & 0 \\
\vdots & \vdots & \ddots & \vdots \\
0 & 0 & \cdots &  2\alpha +1  \\
\end{array}
\end{bmatrix}\mathstrut
\end{array}.$$

Now,

 $$\Delta(\mathcal{Y}, x) =
\begin{vmatrix}

\mathcal{A}_{4 \times 4} - xI_{4 \times 4} & \mathcal{B}_{4 \times 2^{k-2}p}  \\
\mathcal{C}_{2^{k-2}p \times 4} & \mathcal{D}_{2^{k-2}p \times 2^{k-2}p} -  xI_{2^{k-2}p \times 2^{k-2}p}
\end{vmatrix}.$$

Note that, $$ (\mathcal{D}_{2^{k-2}p \times 2^{k-2}p} - xI_{2^{k-2}p \times 2^{k-2}p})^{-1} = \begin{array}{@{} c @{}}
\begin{bmatrix}
\begin{array}{ *{8}{c} }
\frac{1}{2\alpha+1-x} & 0 & \cdots & 0 \\
0 & \frac{1}{2\alpha+1-x} & \cdots & 0 \\
\vdots & \vdots & \ddots & \vdots \\
0 & 0 & \cdots & \frac{1}{2\alpha+1-x}  \\
\end{array}
\end{bmatrix}\mathstrut
\end{array},$$  and $$\begin{vmatrix}
\mathcal{D}_{2^{k-2}p \times 2^{k-2}p} - xI_{2^{k-2}p \times 2^{k-2}p} 
\end{vmatrix} = (2\alpha+1-x)^{2^{k-2}p}.$$ 
So, by Theorem \ref{t6.1}, we have $\Delta(\mathcal{Y}, x)=$
$$ \begin{vmatrix}
 \mathcal{D}_{2^{k-2}p \times 2^{k-2}p} - xI_{2^{k-2}p \times 2^{k-2}p} 
\end{vmatrix}\begin{vmatrix}

\mathcal{A}_{4 \times 4}-xI_{4 \times 4}-\mathcal{B}_{4 \times 2^{k-2}p}(\mathcal{D}_{2^{k-2}p \times 2^{k-2}p} - xI_{2^{k-2}p \times 2^{k-2}p})^{-1}\mathcal{C}_{2^{k-2}p \times 4}
\end{vmatrix}.$$ 

\bigskip
Now,
$|
\mathcal{A}_{4 \times 4}-xI_{4 \times 4}-\mathcal{B}_{4 \times 2^{k-2}p}(\mathcal{D}_{2^{k-2}p \times 2^{k-2}p} - xI_{2^{k-2}p \times 2^{k-2}p})^{-1}\mathcal{C}_{2^{k-2}p \times 4}| = $

$$ \begin{array}{@{} c @{}}
\begin{vmatrix}
\begin{array}{ *{8}{c} } 	
\alpha - x  & 1 - \alpha  & 0 & 0 \\
2^{k-1}p(1-\alpha) & 2^{k + 1}p\alpha - \alpha -  \frac{2^{k-1}p (1 -  \alpha)^2}{2\alpha+1-x} - x
 & (2^kp-2)(1-\alpha) & 1-\alpha - \frac{2^{k - 1} p(1 -  \alpha)^2}{2\alpha+1-x}
  \\
0 & 1-\alpha & 2^{k}p+2\alpha-3-x & 1-\alpha \\
0 & 1-\alpha - \frac{2^{k-1}p (1-\alpha)^2}{2\alpha+1-x} & (2^{k}p-2)(1-\alpha) & 3 \cdot 2^{k-1}p\alpha - \alpha - \frac{2^{k-1}p (1 - \alpha)^2}{2\alpha+1-x}-x \\
\end{array}
\end{vmatrix}\mathstrut
\end{array}$$

$=\frac{1}{2\alpha+1-x} \cdot P(x)$, where 

\bigskip

$P(x) = \bigg[x^5 - \Big[\big(7 \cdot 2^{k-1}p+3\big)\alpha+2^{k}p-2\Big]x^4 - \Big[\big(-3 \cdot 2^{2k}p^2-21 \cdot 2^{k-1}p-2\big)\alpha^2+\big(-7 \cdot 2^{2k-1} p^2-2^{k}p+8\big)\alpha+5 \cdot 2^{k-1}p\Big]x^3- \Big[\big(9 \cdot 2^{2k}p^2+7 \cdot 2^{k}p\big)\alpha^3+\big(3 \cdot 2^{3 k}p^3+23 \cdot 2^{2k-1}p^2-3 \cdot 2^{k+1}p-6\big)\alpha^2+ \big(2^{2k-1}p^2-23 \cdot 2^{k}p+6\big)\alpha-3 \cdot 2^{2k-1}p^2+2^{k+2}p+2\Big]x^2-
\Big[-3 \cdot 2^{2k+1}p^2\alpha^4+\big(-13 \cdot 2^{3k-1}p^3-59 \cdot 2^{2k-2} p^2+5 \cdot 2^{k+1}p\big)\alpha^3+\big(-2^{3k+3}p^3+105 \cdot 2^{2k-2}p^2+45 \cdot 2^{k-1} p-6\big)\alpha^2+ \big(5 \cdot 2^{3k-1}p^3+9 \cdot 2^{2k-2}p^2-5 \cdot 2^{k+2}p\big)\alpha-5 \cdot 2^{2k-2}p^2+3 \cdot 2^{k-1}p+1\Big]x-\big(3 \cdot 2^{3k}p^3+15 \cdot 2^{2k-1}p^2\big)\alpha^4- \big(31 \cdot 2^{3k-2}p^3-95 \cdot 2^{2 k-2}p^2-2^{k+1}p\big)\alpha^3-\big(-2^{3k-2}p^3-37 \cdot 2^{2k-2}p^2+9 \cdot 2^{k+1}p-2\big)\alpha^2-\big(-7 \cdot  2^{3k-2}p^3+25 \cdot 2^{2k-2}p^2-3 \cdot 2^{k-1}p-1\big)\alpha-2^{3k-2}p^3+2^{2k-2}p^2+2^{k}p\bigg].$

Thus from the above, $2\alpha+1-x$ is an $A_{\alpha}$ eigenvalue with $m(2\alpha+1-x) = 2^{k-2}p-1$ and the remaining five eigenvalues are the roots of the polynomial $P(x)$. This completes the proof. \hfill $\Box$

\section{$RD_{\alpha}$ eigenvalues of the power graph $P(\mathcal{G})$}\label{s4}

In this section, we determine the $RD_{\alpha}$ eigenvalues of the power graph $P(\mathcal{G})$.
For this, consider an $n \times n$ symmetric matrix of the form 

\begin{equation} \label{M1}
  \begin{bmatrix}
	\mathcal{U} & \mathcal{V} & \mathcal{V} & \cdots & \mathcal{V} \\
	\mathcal{V}^T & \mathcal{X} & \mathcal{W} & \cdots & \mathcal{W} \\
	\mathcal{V}^T & \mathcal{W} & \mathcal{X} & \cdots & \mathcal{W} \\
	\vdots & \vdots & \vdots &\ddots & \vdots \\
	\mathcal{V}^T & \mathcal{W} & \mathcal{W} & \cdots & \mathcal{X}
	\end{bmatrix}
\end{equation}

where  $\mathcal{U} \in R^{m_1 \times m_1}, \mathcal{V} \in R^{m_1 \times m_2}$,  $\mathcal{X}, \mathcal{W} \in R^{m_2 \times m_2}$, 
and $n = m_1 + cm_2$, $c$ is the number of copies of $\mathcal{X}$.  Then 
the following theorem is used to determine the spectrum of the matrix $RD_{\alpha}$: 

\begin{thm} \label{4.1} \cite{ft(2016)}
	Let the matrix $\mathcal{M}$ is of the same form as given in Equation \ref{M1} and $Spec(\mathcal{M})$ denotes the spectrum of $\mathcal{M}$. Then $Spec(\mathcal{M}) = Spec(\mathcal{N}) \cup Spec(\mathcal{X}-\mathcal{W})^{c-1}$, where $\mathcal{N} = \begin{bmatrix}
	\mathcal{U} & \sqrt{c}\ \mathcal{V}  \\ 
	\sqrt{c}\ \mathcal{V}^T & \mathcal{X} + (c-1)\mathcal{W} & \\ \end{bmatrix}$ and $Spec(\mathcal{X}-\mathcal{W})^{c-1}$ denotes the collection of eigenvalues of $\mathcal{X}-\mathcal{W}$ each having multiplicity $c-1$.
\end{thm}

The proof of below lemmas follows by using the Lemma \ref{4.1}.

\begin{thm}
	
The $RD_{\alpha}$-spectrum of the power graph $P(\mathcal{G})$ is \\

$\Biggl\{\bigg(\alpha(1 + 2^kp)\bigg)^{2^{k-2}p-1}, \bigg(2^kp + 2)\alpha - 1\bigg)^{2^{k-2}p-1},\bigg(\alpha \cdot 2^kp -\frac{1-\alpha}{2}\bigg)^{2^{k-1}p-1},\bigg(
(2^kp+1)\alpha -(1-\alpha)\bigg),\bigg(\alpha(3 \cdot 2^{k-1}p-1)-(1-\alpha)\bigg)^{2^kp-3},x_1,x_2,x_3,x_4,x_5 \Bigg\},$ where $x_i \ (1 \leq i \leq5)$ are the  eigenvalues of the matrix $X = \begin{bmatrix}
S & Y  \\
Y^T & \alpha(3 \cdot 2^{k-1}p-1) + (2^kp-3)(1-\alpha)
\end{bmatrix}, 
$ 

\begin{align*}
S & =\begin{bmatrix}
\begin{array}{@{}cccccccc@{}}
\alpha(2^{k+1}p-1) & 1- \alpha & \sqrt{2^{k-1}p}( 1- \alpha) & \sqrt{2^{k-1}p}( 1- \alpha)  \\
1-\alpha & \alpha(7 \cdot 2^{k-2}p-1) & \frac{\sqrt{2^{k-1}p}( 1- \alpha)}{2} & \sqrt{2^{k-1}p} \ (1 - \alpha) \\
\sqrt{2^{k-1}p}( 1- \alpha) & \frac{\sqrt{2^{k-1}p}( 1- \alpha)}{2} & \alpha \cdot 2^kp + \frac{(2^{k-1}p-1)(1-\alpha)}{2} & \frac{\sqrt{2^{2k-2}p^2}( 1- \alpha)}{2}  \\ 
\sqrt{2^{k-1}p}( 1- \alpha) & \sqrt{2^{k-1}p}( 1- \alpha) & \frac{\sqrt{2^{2k-2}p^2}( 1- \alpha)}{2} & \alpha \cdot 2^kp + \frac{(2^{k-1}p-1)(1-\alpha)}{2} 
\end{array}
\end{bmatrix}
\mathstrut, \\ \bigskip Y  &= \begin{bmatrix}
\sqrt{2^kp-2}(1-\alpha) & \sqrt{2^kp-2}(1-\alpha) & \frac{\sqrt{2^{k-1}p(2^kp-2)}(1-\alpha)}{2} & \frac{\sqrt{2^{k-1}p(2^kp-2)}(1-\alpha)}{2}
\end{bmatrix}^T.
\end{align*}

\end{thm}

\section{Distant and detour distance properties} \label{s5}

In this section, we determine various distant and detour distance properties such as $ec_D, rad_D, \linebreak dia_D$, metric dimension, strong metric dimension.  Now, we show the following lemma:

\begin{lem}\label{l5.1}
Let $P(\mathcal{G})$ be the power graph of $\mathcal{G}$. Then

	\begin{enumerate}
		\item 
		\begin{equation*}
		\hspace{-.5cm}	ec_D(a)= 
		\begin{cases}
		2^{k}p+1 & \text{$ \text{if} \ a \in H_0,$}\\	
		2^{k}p+3 & \text{$ \text{if} \ a \in H_1,$} \\
		2^{k}p+2 & \text{$ \text{if} \ a \in H_2,$}\\
		2^{k}p+3 & \text{$ \text{if}\ a \in H_3,$}
		
		\end{cases}       
		\end{equation*}
		
		\item $rad_D(P(\mathcal{G})) = 2^{k}p+1,$
		
		\item $dia_D(P(\mathcal{G})) = 2^{k}p +3.$
		
	\end{enumerate}
\end{lem}

\noindent {\bf Proof.} For any two distinct vertices of $P(\mathcal{G})$, we have the following possibilities:
\begin{itemize}
	
	\item $e$ is adjacent to all the vertices of $P(\mathcal{G}) \setminus \{e\},$
	\item two distinct vertices of $H_1 \cup \{u\}$ are adjacent while two distinct vertices of $H_2$ are non-adjacent,
	\item $h_1 \in H_1$ is non-adjacent with $h_2 \in H_2,$
	\item $h_1 \in H_1 $ is non-adjacent with $h_3 \in H_3,$
	\item $u \in H_0$ is adjacent with $h_3 \in H_3,$
	\item $h_2 \in H_2$ is non-adjacent with $h_3 \in H_3,$
	\item $sr^{2j+1}$, $sr^{2j+1+2^{k-1}p} \in H_3$  are adjacent for $ 0 \leq j \leq 2^{k-2}p-1$ and $sr^{2j+1},sr^{2i+1} \in H_3$ are non-adjacent for $j \neq i$, $i \neq j + 2^{k-2}p$, $ 0 \leq j \leq 2^{k-1}p-1$.
\end{itemize}

\smallskip
Thus from the above cases, we have

\smallskip

If $a=e$, then $d_D(a,h_1) = 2^kp+1$ for every $h_1 \in H_1$, $d_D(a,u) = 2^kp-1$, $d_D(a,h_2) = 1$ for every $h_2 \in H_2$, and $d_D(a,h_3) = 2^kp+1$ for every $h_3 \in H_3.$ 
Thus, $ec_D(a=e) = 2^{k}p + 1$. 

If $a = u$, then $d_D(a,h_1) = 2^kp+1$ for every $h_1 \in H_1,$ $d_D(a,h_2) = 2^kp$ for every $h_2 \in H_2$, $d_D(a,h_3) = 2^kp+1$ for every $h_3 \in H_3$, and $d_D(u,e) = 2^kp-1$. 
Thus, $ec_D(a = u) = 2^{k}p + 1$. 

If $a \in H_1$, then $d_D(a,h_1) = 2^{k}p+1$ for every $h_1 \in H_1 \setminus \{a\},$ $d_D(a,h_2) = 2^{k}p+2$ for every $h_2 \in H_2,$ $d_D(a,h_3) = 2^{k}p+3$ for every $h_3 \in H_3$, $d_D(a,e) = 2^kp+1$, and $d_D(a,u) = 2^kp+1$. 
Thus, $ec_D(a) = 2^{k}p + 3,$  if $a \in H_1$.

If $a \in H_2$, then $d_D(a,h_2) = 2$ for every $h_2 \in H_2 \setminus \{a\},$ $d_D(a,h_1) = 2^kp+2$ for every $h_1 \in H_1$, $d_D(a,e) = 1$, $d_D(a,u) = 2^kp$, and $d_D(a,h_3) = 2^kp+2$ for every $h_3 \in H_3$. 
Thus, $ec_D(a) = 2^{k}p + 2,$  if $a \in H_2$. 

If $a \in H_3$, then $d_D(sr^{2j+1},sr^{2j+1+2^{k-1}p}) =  2^kp+1$ for $ 0 \leq j \leq 2^{k-2}p-1$ and $d_D(sr^{2j+1},sr^{2i+1}) =  2^kp+3$ for $j \neq i$, $i \neq j + 2^{k-2}p$, where $ 0 \leq j \leq 2^{k-1}p-1$, $d_D(a,u) =  2^kp+1$, $d_D(a,h_1) =  2^kp+3$ for every $h_1 \in H_1$, $d_D(a,e) = 2^kp+1$, and $d_D(a,h_2) =  2^kp+2$ for every $h_2 \in H_2$. 
Thus, $ec_D(a) = 2^{k}p + 3,$  if $a \in H_3$. 

\smallskip 
From the above computation, it follows that the minimum and maximum $ec_D(a)$ of $P(\mathcal{G})$ are $2^kp+1$ and $2^kp+3$, respectively. Therefore, $rad_D(P(\mathcal{G}))= 2^kp+1$ and $dia_D(P(\mathcal{G})) = 2^kp+3$. \hfill $\Box$

\begin{lem}\label{l5.2}
The metric dimension of $P(\mathcal{G})$ is $7 \cdot 2^{k-2}p-4.$
\end{lem}

\noindent{Proof.} Let $S$ be a resolving set of $P(\mathcal{G})$. 
Then, the sets $$\{r,r^2, \ldots, r^{2^{k-1}p-1},r^{2^{k-1}p+1}, \ldots, r^{2^{k}p-1}\}, ~\{sr, sr^{2^{k-1}p+1}\},$$ $$\{sr^3, sr^{2^{k-1}p+3}\}, \ldots, \{sr^{2^{k-1}p-1},sr^{2^{k}p-1}\}, \text{ and } \{s,sr^2,sr^4, \ldots, sr^{2^{k}p-2}\}$$ are twin sets of $P(\mathcal{G})$.
Note that, any resolving set of $P(\mathcal{G})$ can omit at most one element from the twin sets. 
So, $S$ must contains at least $(2^kp-2)-1$, $\underbrace{2-1, 2-1, \ldots, 2-1, 2-1}_{2^{k-2}p-times}$, and $2^{k-1}p-1$ elements of $P(\mathcal{G})$. Thus $\psi(P(\mathcal{G})) \geq 7 \cdot 2^{k-2}p-4
$. 
Moreover, 
$$\left\{ r, r^2, \ldots, r^{2^{k-1}p-1},r^{2^{k-1}p+1}, \ldots, r^{2^{k}p-2},sr^2,sr^4, \ldots, sr^{2^kp-2},sr, sr^3, \ldots, sr^{2^{k-1}p-1} \right\}$$ is a resolving set of order $7 \cdot 2^{k-2}p-4$.
So, $\psi(P(\mathcal{G})) \leq 7 \cdot 2^{k-2}p-4$. 
Hence,$$\psi(P(\mathcal{G})) = 7 \cdot 2^{k-2}p-4.$$ \hfill $\Box$

Now, we calculate the strong metric dimension of $P(\mathcal{G})$, using the following lemma.

\begin{lem} [\cite{of(2007)}] \label{l5.3}
 For any connected graph $G$, $sdim(G) = \beta(G_{SR})$.
\end{lem}

\begin{lem}
		Let $P(\mathcal{G})$ be the power graph of $\mathcal{G}$. Then
		$sdim(P(\mathcal{G})) = 2^{k+1}p-3.$
	\end{lem}
	\textbf{Proof.} Note that $N(x) = (H_0 \ \cup \ H_1)  \setminus \{x\}$ for $x \in H_1$,
	$N(y) = \{e\}$ for $y \in H_2$, $N(e) = \mathcal{G} \setminus \{e\} $, $N(r^{2^{k-1}p}) = (H_0 \ \cup H_1 \ \cup H_3) \setminus \{r^{2^{k-1}p}\}$, $N(z) = \{e, r^{2^{k-1}p}, sr^{2j+1},sr^{2j+1+2^{k-1}p} \} \setminus \{z\}$ for $z \in \{sr^{2j+1},sr^{2j+1+2^{k-1}p}\}$, where $0 \leq j \leq 2^{k-2}p-1.$ Now, we have the following:
	
	\begin{enumerate}
	\item For $u, v \in P(\mathcal{G}) \setminus \{e, r^{2^{k-1}p}\}$ and $w \in N(u)$, $d(v,w) \leq d(u,v)$. Also, for $u, v \in P(\mathcal{G}) \setminus \{e, r^{2^{k-1}p}\}$ and $w \in N(v)$, $d(u,w) \leq d(v,u)$. Thus $u \ MMD \ v$ if $u, v \in P(\mathcal{G}) \setminus \{e, r^{2^{k-1}p}\}$.
	
	\item For $u = r^{2^{k-1}p}$, $v \in H_2$, and $w \in N(r^{2^{k-1}p})$, $d(v,w) \leq d(r^{2^{k-1}p},v)$. Also, for $u = r^{2^{k-1}p}$, $v \in H_2$ and $w \in N(v)$, $d(r^{2^{k-1}p},w) \leq d(v,r^{2^{k-1}p})$. Thus $r^{2^{k-1}p} \ MMD \ v$ if $u = r^{2^{k-1}p}$ and $v \in H_2$.
	
	\item $e$ is not MD with $u \in P(\mathcal{G}) \setminus \{e\}$ as $d(u,w) \nleq d(e,u)$ for every $w \in N(e).$
	
	\item $ r^{2^{k-1}p}$ is not MD with $u \in P(\mathcal{G}) \setminus \{e\}$ as $d(u,w) \nleq d(r^{2^{k-1}p},u)$ for every $w \in N(r^{2^{k-1}p}).$

	\end{enumerate}
	
	This gives that $\mathcal{G}_{SR} \cong K_{2^{k+1}p-2} \cup \{e\} \cup K_{1,2^{k-1}p}$, where $K_{1,2^{k-1}p} \cong {r^{2^{k-1}p} \vee H_3}$. Clearly, $\beta(\mathcal{G}_{SR}) = (2^{k+1}p-2)-1$ as $\beta(K_n) = n-1$. Thus, by Lemma \ref{l5.3}, $sdim(P(\mathcal{G})) = \beta(\mathcal{G}_{SR}) = 2^{k+1}p-3.$  \hfill $\Box$

\begin{lem}\label{l5.4}
	Let $P(\mathcal{G})$ be the power graph of $\mathcal{G}$.
	Then
	
	\begin{enumerate}
		\item 
		$dds(\mathcal{G}) = \Big(\big(1, 2^{k+1}p-1 \big), \big( 1,3 \cdot 2^{k-1}p-1,2^{k-1}p \big), \big( 1, 2^{k}p -1, 2^{k}p \big)^{2^kp-2}\Big),$ 
		
		\item  $dds_D(G) = \Big( \big(1,2^{k-1}p,0^{2^{k}p-3},1,0,3 \cdot 2^{k-1}p-2 \big), \big( 1,0^{2^{k}p-2},1,2^{k-1}p,3 \cdot 2^{k-1}p-2 \big), \linebreak \big( 1,0^{2^{k}p},2^{k}p-1, (2^{k-1}p)^2 \big)^{2^{k}p-2}, \big(1^2,2^{k-1}p-1,0^{2^{k}p-3},1,0,3 \cdot 2^{k-1	}p-2 \big)^{2^{k-1}p},\big(1,0^{2^{k}p},3, \\ 2^{k-1}p,3 \cdot 2^{k-1}p-4 \big)^{2^{k-1}p}\Big)$.
		
	\end{enumerate}
\end{lem}

\noindent {\bf Proof.} For $a \in V(P(\mathcal{G}))$, we get 

\begin{equation*}
ec(a) = \begin{cases}
1 & \text{$ \text{if} \ a = e,$}\\
2 & \text{$ \text{if} \ a \neq e.$}\\
\end{cases}
\end{equation*}

Also, by Lemma \ref{l5.1}, we have

\begin{equation*}
\hspace{-.5cm}	ec_D(a)= 
\begin{cases}
2^{k}p+1 & \text{$ \text{if} \ a \in H_0,$}\\	
2^{k}p+3 & \text{$ \text{if} \ a \in H_1,$} \\
2^{k}p+2 & \text{$ \text{if} \ a \in H_2,$}\\
2^{k}p+3 & \text{$ \text{if}\ a \in H_3.$}

\end{cases}       
\end{equation*}

Therefore, \begin{equation*}
dds(a)= 
\begin{cases}
(1, 2^{k+1}p-1) & \text{$ \text{if} \ a = e,$}\\
(1,3 \cdot 2^{k-1}p-1,2^{k-1}p) & \text{$ \text{if} \ a = u,$}\\
(1, 2^{k}p -1, 2^{k}p) & \text{$ \text{if} \ a \in H_1$.}

\end{cases}       
\end{equation*}

\begin{equation*}
dds_D(a)= 
\begin{cases}
(1,2^{k-1}p,0^{2^{k}p-3},1,0,3 \cdot 2^{k-1}p-2) & \text{$ \text{if} \ a = e,$}\\
(1,0^{2^{k}p-2},1,2^{k-1}p,3 \cdot 2^{k-1}p-2) & \text{$ \text{if} \ a = u$,}\\
(1,0^{2^{k}p},2^{k}p-1, (2^{k-1}p)^2) & \text{$ \text{if} \ a \in H_1$,} \\
(1^2,2^{k-1}p-1,0^{2^{k}p-3},1,0,3 \cdot 2^{k-1}p-2) & \text{$ \text{if} \ a \in H_2$,}\\
(1,0^{2^{k}p},3,2^{k-1}p,3 \cdot 2^{k-1}p-4) & \text{$ \text{if} \ a \in H_3$.}

\end{cases} 
\end{equation*}

Using the fact that,  $|\{e\}|$, $|\{u\}|,$  $|H_1|$, $|H_2|$, and $|H_3|$ are  $1,1, 2^{k}p-2, 2^{k-1}p, \text{ and } \ 2^{k-1}p$, respectively, we get the required proof. \hfill $\Box$

\section{Acknowledgment}

\end{document}